\documentclass[12pt]{article}
\usepackage{amssymb}
\newtheorem{thm}{ Theorem}
\newtheorem{defn}{ Definition}
\newtheorem{lem}{ Lemma}

\newtheorem{rem}{\it Remark}[section]

\begin{document}
\baselineskip=22pt
\title{\bf Almost  periodic  solutions of differential equations with piecewise constant argument of generalized type}
\author{M. U. Akhmet\thanks{M.U. Akhmet is previously known as M. U. Akhmetov.}} 
\date{{\small Department of Mathematics and Institute of Applied Mathematics, Middle East 
Technical University, 06531 Ankara, Turkey}\\{\bf marat@metu.edu.tr}}

\maketitle

\vspace{0.5cm}

\
\begin{abstract}
We consider existence and stability  of an
almost periodic solution of the following hybrid system 
\begin{eqnarray}
\frac{dx(t)}{dt} = A(t)x(t) + f(t,x(\theta_{\beta(t) - p_1}) ,x(\theta_{\beta(t) - p_2}), \ldots, x(\theta_{\beta(t) - p_m})),
\label{e1}
\end{eqnarray}
where $x \in \mathbb R^n, t \in \mathbb R, \beta(t) = i$ if $\theta_i \leq t <\theta_{i+1}, i= \ldots -2, -1, 0, 1, 2, \ldots,$ is an identification function, 
$\theta_i$ is a strictly ordered sequence of real numbers, unbounded on the left and on the right, $p_j, j=1,2, \ldots,m,$ are fixed integers,  and the 
linear homogeneous system associated with (\ref{e1}) satisfies exponential dichotomy. The  deviations of the argument are not restricted by any sign assumption when existence is considered. The problem  of positive (almost periodic) solutions of the logistic equation  is discussed as an example. 
 A new technique  of investigation of equations with piecewise argument, based on integral representation, is developed. 
\end{abstract}

\noindent {\it Key words and phrases:} Quasilinear system; Almost 
periodic solutions;  Piecewise constant argument of general type; Advanced-delayed 
argument; Differential logistic equations; Positive solutions.

\noindent {\it 1991 Mathematics Subject Classification:} 34K14; 34K20.

\vspace{0.5cm}
\section{Introduction and Preliminaries}
The theory of differential equations with  piecewise constant argument  (EPCA) of the type
\begin{eqnarray}
&& \frac{dx(t)}{dt} = f(t,x([ t - p_1]) ,x([ t - p_2]),\ldots, x([ t - p_m])),
\label{e2}
\end {eqnarray}
where $[\cdot]$ signifies the greatest integer function, was initiated in  \cite{cw1, shah} and developed by many authors
 \cite{aw1, aa, aa1, ky, p, p1}, \cite{s}, \cite{w2}-\cite{ym}. These systems have been under intensive investigation for the last twenty years. They describe hybrid dynamical systems
and combine properties of both differential and difference equations. 
An example of  the application of these equations to the problems of  biology can be found in \cite{bus}.
One of the novel ideas in our paper is that  system (\ref{e1})  is of general type (EPCAG)  for equation (\ref{e2}). 
Indeed if we take $\theta_i = i, i = \ldots, -2, -1, 0, 1, 2, \ldots $   then   (\ref{e1}) takes the form of    (\ref{e2}).

 The existing method of investigation  of EPCA, as  it was proposed 
 by its founders \cite{cw1, shah}, is based on the reduction of EPCA to discrete equations, and  it has been the only method to prove assertions about EPCA until now. In our paper, we (apparently for the first time)  propose another approach to the problem. In fact,  we are 
 dealing with the construction  of an equivalent integral equation. Since we do not need additional assumptions on the reduced discrete equations for investigating EPCAG, the new method requires more easily verifiable conditions,  similar to those for ordinary differential  equations. So, solving the problems of EPCAG (as well as of  EPCA) may become less cumbersome if the approach proposed in our paper is applied. 
 
 Another novelty in our investigation is that we consider  equations with deviated argument
of mixed (advanced-delayed) type.  Even in the case of  advanced argument,  there are certain difficulties  if we try to define a solution for increasing  $t$  \cite{e}. J. Hale remarked in \cite{h2} that \rq\rq these
equations (of mixed type) seem to dictate that boundary conditions should be specified in order to obtain a solution
in the way as one does for elliptic partial differential equations.\lq\lq  \, We  regard the boundedness of the solution on $R$ as a  boundary condition in our investigation. Similar arguments were used in \cite{ap6,b,dn} to investigate various problems for ordinary and functional differential equations. 

 The existence of almost periodic solutions is one of the most 
interesting subjects  of the theory of differential equations 
(see, for example, \cite{amerio, cord1, f} and the references cited there).
This problem has been  considered in the context of EPCA in many papers, such as \cite{s, y1, z}.

To solve the problems of the present paper,  we intend to apply our knowledge about the almost periodicity of discontinuous solutions of impulsive systems
\cite{ap}-\cite{aps}, \cite{hw,sp,slyus}. 
One should not be surprised 
with the relation between EPCAG and impulsive differential equations.
This possibility  was mentioned in \cite{cw1} for EPCA,  and in \cite{fil} for differential equations with discontinuous right hand side.

Original ideas on the spaces of discontinuous functions are to be found in \cite{hw,kol,skor,w}. Following these results,  in \cite{ap1} we introduced  Bohner type discontinuous almost periodic 
functions using topology as well as metric, in the spaces of discontinuous functions and  in the discrete spaces of sets of points on the real axes, unbounded on the left and on the right. The multiplicity of one for the elements of the sets was mentioned there, since we supposed that the distances between neighbors are uniformly separated from zero. Further our  proposals  based on  the metric were developed in  Supplement $A$ of  \cite{sp}.  In the present paper we again consider  the spaces from the topological point of view, assuming that the multiplicity more or equal to one.

One can be confident  that the reduction to integral equations, as well as the  awareness about  the theory of discontinuous functions, can diminish the number of\\ \lq\lq {strange} \rq\rq 
properties of EPCA, which  are  usually generated by the reduction to difference equations,  and can give  explanations of certain phenomena. For example,  the result on the module containment considered in \cite{y1} becomes less specific if one compares it with our Theorem \ref{thm1} (see also Example 1). But we should note that the reduction to discrete equations is preferable in some cases, as in \cite{g,K,m}, where the period-doubling bifurcation and the generation of chaos by a logistic EPCA are considered.

Let $ \mathbb  Z,\mathbb N,$ and $\mathbb R$ be the sets of all integers, natural and real numbers, respectively, and
$||\cdot||$ be the euclidean norm in $\mathbb R^n, n \in \mathbb N.$ 
 Let   $s\in R$ be a positive number. We denote  $ G_s = \{x \in \mathbb R^n | ||x|| \leq s\}$  and 
$G_s^{m} = G_s \times G_s \ldots  \times G_s$ (that is, $G_s^{m}$ is an $m-$ times Cartesian 
product of $G_s).$
Let a  $C_0(R)$ (respectively $C_0(R \times G_H^{m})$ for a given $H \in R,\,H >0)$ be the set of all bounded and uniformly continuous functions on $R$ (respectively on $R\times G_H^{m}).$
For $f\in C_0(R)$ (respectively $C_0(R\times G_H^{m})$) and $ \tau \in R,$ 
the translate of $f$ by $\tau$ is  the function
$Q_{\tau} f = f(t+\tau), \, t \in R$ ( respectively $Q_{\tau} f(t,z)  = f(t+\tau,z),\, (t,z)  \in R\times G_H^{m}).$ 
 A number  $\tau \in R$ is called an $\epsilon -$ translation number of a function
 $f\in C_0(R)$ ( $C_0(R \times G_H^{m})$) if  $\,|| Q_{\tau} f - f||< \epsilon\,$ for every
 $t \in R \, (\, (t, z) \in R \times G_H^{m}).$ A set $S \subset R$  is said to be relatively dense if there 
exists a number $l > 0,$ such that  $[a, a+ l] \cap S \not = \emptyset $ for all $a \in R.$ 
\begin{defn}\rm   A function $f \in C_0 (R) (C_0(R \times G_H^{m}))$ is said to be almost periodic 
 (almost periodic in $t$ uniformly with respect to 
 $ z \in G_H^{m}$)  if for every  $\epsilon \in R,  \epsilon > 0,$  there exists  a relatively dense set of 
$\epsilon -$ translation numbers of f.
 \end{defn}
     Denote   by ${\cal AP}(R)$ (${\cal AP}(R\times G_H^{m})$) the set of all such functions .\\
The following assumptions will be needed throughout the paper.
\begin{itemize}
\item[($C_1$)] $A(t)\in {\cal AP}(R)$ is an $n \times n$ matrix;
\item[($C_2$)] $f\in {\cal AP}(R\times G_s^{m}),$ for every $s\in R, s \geq 0;$
\item[($C_3$)] $\exists \,l\in R, l>0,$ such that 
$$ ||f(t,z_1) - f(t,z_2)|| \leq l \sum_{j=1}^{m} ||z_1^j - z_2^j||,$$
 where  $z_i = (z_i^1,\ldots, z_i^{m}) \in R^{nm}, i=1,2.$ \\
Let 
\begin{eqnarray}
&& \frac{dx}{dt} =  A(t)x 
\label{e3}
\end{eqnarray}
be the  homogeneous 
linear system associated with (\ref{e1}),  and  $X(t)$ be a fundamental matrix of (\ref{e3}).
\item[($C_4$)] system  (\ref{e3}) satisfies  exponential dichotomy, that is,  there exist
a projection $P$ and positive constants $\sigma_1, \sigma_2, K_1, K_2,$ such that
\begin{eqnarray*}
&&||X(t)PX^{-1}(s)|| \leq K_1 \exp(- \sigma_1 (t-s)), \, t \geq s,\nonumber\\
&&||X(t)(I-P)X^{-1}(s)|| \leq K_2 \exp( \sigma_2 (s-t)), \, t \leq s.
\end{eqnarray*}
\end{itemize}

Let $\Theta$ be  a space of  strictly ordered sequences $\{\theta_i\} \subset R, i \in \mathbb Z,$ such that
$|\theta_i| \rightarrow \infty,$ if  $|i| \rightarrow \infty.$  Denote by 
${\cal PC}$ the set of all functions from $\mathbb R$ to  $\mathbb R^n$ that are piecewise continuous with discontinuities of the first type.
Assume  that the set of discontinuities of every function from ${\cal PC},$ numerated in a strict order, is  an element 
of $\Theta.$ Moreover, these functions are uniformly continuous on the set $\cup_{i \in \mathbb Z}(\theta_i,\theta_{i+1}),$ and they are left or right continuous at every point of discontinuity.\\
Denote by ${\cal PC}_r \subset {\cal PC}$ the set of all continuous from the right functions. Similarly, one can define a set
${\cal PC}_l.$ If $\phi \in {\cal PC}$, then  one can  define a function $\phi_r \in {\cal PC}_r,$  such that
$\phi_r (t) = \phi(t)$ everywhere, except possibly at points $t = \theta_i,$  that is,

$$ \phi_r(t) =\left\{\begin{array}{ll}
\phi(t), & \mbox {if $t \not = \theta_i$},\\
\phi(\theta_i+), i \in Z , & \mbox{otherwise.}
\end{array}\right.$$
We shall call the function $\phi_r(t)$  a right extension of the function $\phi(t) \in {\cal PC}.$  Since the function  $\beta(t)$ is right continuous,  it is reasonable to consider only the space ${\cal PC}_r$  in our paper,  extending, to the right,  if necessary, every function from ${\cal PC}$ that we obtain in our discussion.
Since functions from ${\cal PC}$ are assumed to be the derivatives or limits of the solutions of EPCAG, no difficulty arises from this agreement.
In what follows we assume that $\beta(t) \in {\cal PC}_r.$
 The following definition of a solution of EPCAG which is a slightly changed form of the corresponding definition  for EPCA  \cite{p,p1} can be given.
\begin{defn}   A function $x(t)$ is a solution of (\ref{e1}) on $\mathbb R$ if:
\begin{enumerate}
\item[(i)] $x(t)$ is continuous on $\mathbb R;$
\item[(ii)] the derivative $x'(t)$ exists at each point $t \in \mathbb R,$ with the possible exception of the points $\theta_i , i \in \mathbb Z,$  where one-sided derivatives exist;
\item[(iii)] equation (\ref{e1}) is satisfied on each interval $[\theta_i, \theta_{i+1}), i \in \mathbb Z.$ 
\end{enumerate}
\label{defn5}
\end{defn}
It is obvious that the derivative of  a solution $x(t)$ is a function from ${\cal PC}_r,$
if we assume  it to be the right derivative at $t= \theta_i, i \in \mathbb Z.$

Let 
$$G(t,s) =\left\{\begin{array}{rr}
 X(t)PX^{-1}(s),  \quad \mbox {if} \quad t \geq s,\\
 X(t)(P-I)X^{-1}(s)  ,\quad \mbox{if} \quad t < s
\end{array}\right. $$
be the Green's function of (\ref{e3}).
Denote
$$F_{\theta}(\psi(t)) = f(t,\psi(\theta_{\beta(t) - p_1}) ,\psi(\theta_{\beta(t) - p_2}), \ldots, 
\psi(\theta_{\beta(t) - p_m})),$$ where $\psi(t) \in C_0(R).$

The following is one of the most important assertion  for our method of 
investigation of EPCAG.
\begin{lem} A function $x(t) \in C_0(R)$ is a solution of (\ref{e1})  if and only if
\begin{eqnarray}
&& x(t)  = \int_{-\infty}^{\infty} G(t,s) F_{\theta}(x(s))ds.
\label{e4}
\end {eqnarray}
\label{l1}
\end{lem}
{\it Proof.} 

\rm {\it Necessity.} Assume that $x(t) \in C_0(R)$ is a solution of (\ref{e1}). Denote

\begin{eqnarray}
&& \phi(t)  = \int_{-\infty}^{\infty} G(t,s) F_{\theta}(x(s))ds.
\label{e5}
\end {eqnarray}
By straightforward calculation we can see that the function $\phi(t)$ is bounded  and  continuous on $R.$

Assume that $t \not = \theta_i, i \in Z.$ Then

\begin{eqnarray*}
&& \phi'(t)  = A(t)\phi(t) +  F_{\theta}(x(t))
\end {eqnarray*}
and 
\begin{eqnarray*}
&& x'(t)  = A(t)x(t) +  F_{\theta}(x(t)).
\end {eqnarray*}
Hence, 
\begin{eqnarray*}
&& [ \phi(t) -x(t)]'  = A(t)[ \phi(t) -x(t)]  .
\end {eqnarray*}
Calculating the limit values at  $t=\theta_j, j \in \mathbb Z,$ we find that
$$\phi'(\theta_j \pm 0)  = A(\theta_j \pm 0)\phi(\theta_j \pm 0) +  F_{\theta}(x(\theta_j \pm 0))$$
$$x'(\theta_j \pm 0)  = A(\theta_j \pm 0)x(\theta_j \pm 0) +  F_{\theta}(x(\theta_j \pm 0)).$$
Consequently,
$$ [ \phi(t) -x(t)]'|_{ t=\theta_j+0} =    [ \phi(t) -x(t)]'|_{ t=\theta_j-0}.$$

Thus,  $[ \phi(t) -x(t)]$ is a continuously differentiable function on $R,$  satisfying 
(\ref{e3}). That is,  $[\phi(t) -x(t)] = 0$ on $R.$\\
{\it Sufficiency.}  Suppose that (\ref{e4}) is valid and $x(t) \in C_0(R).$ Fix $i \in \mathbb Z$ and consider the
interval $[\theta_i, \theta_{i+1}).$ If  $ t \in (\theta_i,\theta_{i+1}),$ then by differentiating one can see that $x(t)$ satisfies (\ref{e1}).
Moreover, considering $ t \rightarrow \theta_i+,$  and taking into account that $\beta(t)$
is a right-continuous function, we obtain that $x(t)$ satisfies (\ref{e1}) on $[\theta_i, \theta_{ i+1}).$ The lemma is proved.
\section{Wexler sequences}
Fix $\theta \in \Theta,$ and consider a sequence $\gamma_i,
i \in \mathbb Z, \gamma_{i+1} \geq \gamma_i,$ such that  for every  $\gamma_i \in \gamma$ 
there exists an element  $\theta_j \in \theta$ such that $\gamma_i = \theta_j.$ Let  $m(i), i \in \mathbb Z,$ be 
the number of elements of $\gamma$ which are equal to $\theta_i.$ We shall call this number the multiplicity 
of $\theta_i$ with respect to $\gamma.$ Denote $m(\gamma) = \sup_i m(i).$  Denote by $\Gamma$ the set of all
sequences $\gamma$ such that $m(\gamma) <\infty.$ If $\gamma \in \Gamma,$ then we shall say that 
$m(\gamma)$ is the maximal multiplicity of $\gamma.$ It is obvious that $|\gamma_i| \rightarrow \infty$ if 
$|i| \rightarrow \infty,$ for every $\gamma \in \Gamma,$ and  that $\Theta \subset \Gamma.$\
We shall call $\theta$ a support of $\gamma$, and $\gamma$ a representative of $\theta$.
Introduce the following distance 
$||\gamma^{(1)} -  \gamma^{(2)}|| = \sup_i || \gamma^{(1)}_i -  \gamma^{(2)} _i||$ if $\gamma^{(1)}, \gamma^{(2)}  \in \Gamma.$
We shall say that elements $\theta^{(1)},\theta^{(2)} \in \Theta$ are
$\epsilon-$ equivalent and write  $\theta^{(1)} \epsilon \theta^{(2)},$ if there exist the representatives $\gamma^{(1)}$ and $ \gamma^{(2)}$ in  $\Gamma$ of  $\theta^{(1)}$ and $\theta^{(2)},$ respectively,
such that $||\gamma^{(1)} - \gamma^{(2)}|| < \epsilon.$ Moreover,  we shall say that these sequences are
in the $\epsilon-$ neighborhoods of each other.\\
The topology  defined on the basis of all $\epsilon-$ neighborhoods, $0<\epsilon< \infty,$ of
all elements of $\Theta$ is named as $B^s-$ topology. Obviously, it is a Hausdorff  topology.
\begin{lem} If $\quad \theta^{(1)} \epsilon_1 \theta^{(2)}, \theta^{(2)} \epsilon_2 \theta^{(3)},$ then 
 $\quad \theta^{(1)} (\epsilon_1+ \epsilon_2) \theta^{(3)}.$
\label{l2}
\end{lem}
{\it Proof.} Let $\gamma^{(1)},  \gamma^{(3)}$ be the  representatives of  $\theta^{(1)},\theta^{(3)},$
respectively, and $\gamma^{(2)},  \gamma^{(4)}$ be the representatives of $\theta^{(2)}$ such that $||\gamma^{(1)} - \gamma^{(2)}|| < \epsilon_1$  and  $||\gamma^{(4)} - \gamma^{(3)}|| < \epsilon_1.$
Let  $m^{(2)},  m^{(4)}$ be the maximal multiplicities of  $\gamma^{(2)},  \gamma^{(4)},$ respectively.
Denote $m_0 = \max (m^{(2)},  m^{(4)})$  and define a representative  $\bar \gamma^{(2)}$ of $\theta^{(2)}$ with
 multiplicity 
$m(i) = m_0, i \in \mathbb Z.$ Enlarging,  if necessary,  the multiplicity of elements and shifting the indeces of elements by the same
number,  in accordance with the change  from $\gamma^{(2)}$ to $\gamma^{(4)},$  we can  construct representatives
 $\bar \gamma^{(1)}, \bar \gamma^{(3)}$  of  $\theta^{(1)},\theta^{(3)},$ respectively, such that $|\bar \gamma^{(2)}_i - \bar \gamma_ i^{(1)}| < \epsilon_1,  |\bar \gamma^{(2)}_i - \bar \gamma_ i^{(3)}| < \epsilon_2, i \in \mathbb Z.$ Consequently,
$||\bar \gamma^{(3)} - \bar \gamma^{(1)}|| < \epsilon_1 + \epsilon_2.$ The lemma is proved.

Let $a_i, i \in \mathbb Z,$ be a sequence in $\mathbb R^n.$  An integer $p$ is called an $\epsilon-$almost period of the sequence, if $||a_{i+p} - a_i|| < \epsilon$ for any $i \in \mathbb Z.$ 
\begin{defn} A sequence $a_i, i \in \mathbb Z,$ is  almost periodic,  if for any 
$\epsilon >0$ there exists a relatively dense set of its $\epsilon-$almost periods.
\end{defn}
Let $\gamma \in \Gamma, i,j \in \mathbb Z.$ Denote $\gamma_i^j = \gamma_{i+j} -  \gamma_{i}$ and define  sequences
$\gamma^j = \{\gamma_i^j \}\in \Gamma, j \in \mathbb Z.$
\begin{defn} \cite{hw, sp} Sequences $\gamma^j, j \in \mathbb Z,$ are equipotentially almost periodic   if  for an arbitrary $\epsilon >0$ there exists
a relatively dense set of $\epsilon-$ almost periods  that are common for all  $\gamma^j, j \in \mathbb Z$.
\end{defn}
\begin{defn} We shall say that $\theta \in \Theta$ is a Wexler sequence, if there exists a representative $\gamma$ of  $\theta$  with 
 equipotentially almost periodic $\gamma^j, j \in \mathbb Z$.
\end{defn}
Let  $\gamma \in \Gamma, \epsilon >0,$ be given. Denote by $T_{\epsilon} \subset R$  the set of numbers $\tau$,
for which  there exists at least one number $q_{\tau} \in   \mathbb Z,$ such that 
\begin{eqnarray}
&&  |\gamma_i^{q_{\tau}} - \tau| < \epsilon , i \in \mathbb Z.
\label{e9}
\end{eqnarray}
Denote by $Q_{\tau}$ the set of all numbers $q_{\tau}$ satisfying (\ref{e9}) for fixed $\epsilon$ and $\tau$, and 
$Q_{\epsilon} = \bigcup_{\tau \in T_{\epsilon}}Q_{\tau}.$ 
The following lemmas were proved in \cite{hw,w} for $m(\gamma) = 1.$ But one can 
easily, repeating the proof in \cite{hw}, to verify that they valid   
 if  $ 1< m(\gamma) < \infty .$ For example, equivalence of conditions $(a)$ and 
 $(b)$ of Lemma \ref{l3} is considered in   \cite{sp}.
\begin{lem} The following statements are equivalent 
\begin{itemize}
\item[(a)] the sequences $\gamma^j, j\in \mathbb Z,$ are equipotentially almost periodic;
\item[(b)] the set $T_{\epsilon}$ is relatively dense for any $\epsilon >0;$
\item[(c)] the set $Q_{\epsilon}$ is relatively dense for any $\epsilon >0.$
\end{itemize}
\label{l3}
\end{lem}
 
\begin{lem} Assume that  sequences  $\gamma^j, j \in \mathbb Z,$ are equipotentially almost periodic. Then
 for arbitrary $l>0$ there exists $n_0 \in N,$ such that  any interval of length  $l$ contains at most $n_0$
elements of $\gamma.$
\label{l4}
\end{lem}

 Fix  $\theta \in \Theta$ and let $h=\{h_n\}, n \in \mathbb N,$ be a sequence of real numbers. Assume that the sequence 
of shifts $\{\theta + h_n\}_n$ is convergent in $B^s-$ topology. We shall denote  the limit element as $Q_h\theta.$
\begin{defn} An element $\theta \in \Theta$ has the Bohner property, if every sequence $\{h^{'}_n\}$ contains a subsequence
$h \subset h^{'},$ such that   $Q_h\theta$ exists.
\end{defn} 
\begin{thm} If $\theta \in \Theta$ is a Wexler sequence, then it satisfies the Bohner property and 
$Q_h\theta$ is a Wexler sequence for arbitrary $h.$
\label{t1}
\end{thm}
{\it Proof.} Let  $\theta$ be a Wexler sequence and  $\gamma^j, j \in \mathbb Z,$ be equipotentially almost periodic, 
where $\gamma$ is the representative of $\theta.$
The almost periodicity of $\gamma^1$ implies that there exists $\kappa >0,$ such that  $0\leq \theta^1_i < \kappa, 
i \in \mathbb Z.$ Hence, for arbitrary $n \in N$ there exists $i_n$ such that $\gamma_{i_n} +h_n \in [0, \kappa].$
Denote $\gamma^{(n)} = \{ \gamma _{i+i_n}\}_i.$
Clearly,    $\gamma^{(n)} \in \Gamma, m(\gamma^{(n)}) = m(\gamma),$ and  $\gamma^{(n)j}, j \in \mathbb Z, $ are
equipotentially almost periodic.\\
Using Theorem 1, p. 129 \cite{hw} and the equipotentially almost periodicity of  $\gamma^{(n)j}, $
one can show that  there exists a subsequence $n_k$, let us say it is $n$ itself, such that  for arbitrary 
$\epsilon >0$ there exists $n(\epsilon)  \in N,$ such that 
\begin{eqnarray}
&& || \gamma^{(m)j} - \gamma^{(p)j}|| < \frac{\epsilon}{2}, j \in \mathbb Z, \mbox{if}\,  m, p >n(\epsilon).
\label{e12''}
\end{eqnarray}
Moreover, without loss of generality , we assume that 
$\gamma_0^{(n)} + h_n \rightarrow \gamma_0^{(0)} \in [0,\kappa].$ Consequently,
for arbitrary $\epsilon >0$ there exists  $n(\epsilon)$ such that  
\begin{eqnarray}
&& | \gamma^{(m)} + h_m - \gamma^{(p)} - h_p|< | \gamma_0^{(m)} + h_m - \gamma_0^{(p)} - h_p|  \nonumber\\
&&+ | \gamma^{(m)i} - \gamma^{(p)i} | < |   < \frac{\epsilon}{2} + \frac{\epsilon}{2}= \epsilon,
\label{e13}
\end{eqnarray}
if $m, p > n(\epsilon).$ That is, if we fix $i \in \mathbb Z,$ then $\{\gamma_i^{(n)} + h_n\}_n$ is a Cauchy sequence, and hence
$\gamma_i^{(n)} + h_n \rightarrow \gamma_i^{(0)}, i \in \mathbb Z.$ 
Furthermore, by (\ref{e13}) the convergence is uniform in $i$ and $  \gamma_{i+1}^{(0)}\geq \gamma_i^{(0)}, i \in \mathbb Z.$ 
Finally, the condition $m(\gamma^{(n)}) = m(\gamma),$ and Lemma \ref{l4}  imply that 
$m(\gamma^{(0)})\leq n_0m(\gamma)<\infty.$ 
It is  obvious that $|\gamma_i^{(0)}| \rightarrow \infty.$ Hence,  $\gamma^0 \in \Gamma.$ 
Assume that $\theta^0$ is the support of $\gamma^{(0)}$. Since  $\gamma^{(n)j}$ are equipotentially almost periodic
and  $| \gamma_i^{(n)j} - \gamma_i^{(0)}| \rightarrow 0$ uniformly in $i$ as     $n \rightarrow \infty,  j\in \mathbb Z,$ one can show that  
$\gamma^{(0)j}, j\in \mathbb Z,$ are equipotentially almost periodic in a similar manner as in the proof of the theorem on almost periodicity of a limit function \cite{f}.
Consequently, $\theta^{(0)}-$ is a Wexler  sequence. \\The theorem is proved.
\begin{thm} $\theta \in \Theta$ is a Wexler sequence if and only if it satisfies the Bohner property.
\label{t2}
\end{thm}
{\it Proof.} {\it Necessity} is proved by Theorem \ref{t1}.\\
{\it Sufficiency.} Assume that $\theta \in \Theta$ is not a Wexler sequence. Then $\theta^j,  j\in \mathbb Z,$ are not 
equipotentially almost periodic,  and by Lemma \ref{l3} there exists a number $\epsilon_0$ and a sequence of sections 
$I_n = [h_n - l_n, h_n + l_n], n \in \mathbb N,$ where  $l_1$ is arbitrary,  $l_n > \max_{m<n}|h_m|$ and the following
inequality
\begin{eqnarray}
&& \sup_{q,k \in \mathbb Z} |\theta_k^q - \xi| \geq \epsilon_0, \xi \in \cup_{n} I_n,
\label{e14}
\end{eqnarray}
is valid.
Consider a sequence of shifts $\quad\theta+ h_n, n\in \mathbb N,$ and denote $h^{'} = \{h_n\}.$
For arbitrary  $m,p \in N, m>p,$  we have that $ h_m-h_p \in I_m$ and 
$ \sup_{i,j \in \mathbb Z} |\theta_i^j - (h_m-h_p)| \geq \epsilon_0$ or
\begin{eqnarray}
&&  \sup_{i,j \in \mathbb Z} |\theta_i +  h_p - \theta_{i+j} -h_m)| \geq \epsilon_0.
\label{e15''}
\end{eqnarray} 
The last  inequality means that $\theta + h_m$ is not in  the $\epsilon_0-$  neighborhood of $\theta + h_p.$ Assume that there exists a subsequence $h \subset h^{'}$ such that  $ \theta + h_{n_k}$  convergens to $\theta^{(0)} \in \Theta$ uniformly in $B^s-$ topology. Then there exist numbers 
$n_m$ and $n_p, n_m>n_p,$ such that $ (\theta + h_{n_p})\,\frac{\epsilon_0}{2}\, \theta^{(0)} $ and  $ (\theta + h_{n_m})\, \frac{\epsilon_0}{2}\,\theta^{(0)} .$ 
By Lemma \ref{l2},   $ (\theta + h_{n_p})\, \epsilon_0 \, (\theta + h_{n_m}).$ The contradiction proves the theorem.

Subsequences $h$ and $g$ are common subsequences  of  sequences $h^{'}$ and $g^{'}$
\,if $h_n = h_{n(k)}^{'}$ and $g_n = g_{n(k)}^{'}$ for some given function $n(k)$ \cite{f}.
The following theorem is an analogue  of Theorem  1.17 
from  \cite{f}, and we shall follow the proof of the theorem presented there.
\begin{thm} A sequence $\theta \in \Theta$ is a Wexler sequence if and only if  for arbitrary 
$h^{'}$ and $g^{'}$ there exist common subsequences $h \subset h^{'}$ and $g \subset g^{'}$ such that 
\begin{eqnarray}
&& Q_{h+g} \theta = Q_h Q_g \theta.
\label{e15'}
\end{eqnarray}
\label{t3}
\end{thm}
{\it Proof.} \\
{\it Necessity.} Assume that $\theta$ is a Wexler sequence. By the previous theorem
 there exists a subsequence  $g^{''} \subset g^{'}$ such that $Q_{g^{''}}\theta$ exists and the limit
 is  a Wexler sequence. Denote  $\eta = Q_{g^{''}}\theta .$  Moreover, if $h^{''} \subset h^{'}$ is common with  $g^{''}$, then one can find a sequence 
$h^{'''} \subset   h^{''}$ such that $ \mu =  Q_{h^{'''} } \eta $ is a Wexler sequence. If $g^{'''} \subset g^{''}$ is common with 
$h^{'''}$ then there exist common subsequences $h \subset h^{'''},  g \subset g^{'''}$ such that
$Q_{g+h} \theta = \zeta.$  Since    $h \subset h^{''},  g \subset g^{''},\quad Q_{g} \theta = \eta,$ $Q_{h} \eta = \mu.$ 
Thus,
 if  $\epsilon>0$ is fixed, then  for sufficiently large $n$ we have 
$\zeta \epsilon (\theta +h_n+g_n),\quad \eta \epsilon (\theta+g_n),$\quad and $\quad \mu \epsilon (\eta + h_n).$
Using Lemma \ref{l2},  we can conclude that $\zeta 3\epsilon \mu.$ Hence, $\zeta = \mu,$  as $\epsilon$ is arbitrarily
small.  \\
{\it Sufficiency.} Suppose that $h^{'}$ is the  given sequence. Taking $g^{'} = 0$ we see that the condition
$Q_{g+h}= Q_g Q_h \theta$ implies that $Q_h \theta$ exists,  and  hence  by Theorem \ref{t2} the sequence $\theta$
is a Wexler sequence. The theorem is proved.

\section{Bohr-Wexler almost periodic functions}
\begin{defn} Let $u_1, u_2 \in {\cal PC}_r,$ and $\theta^{(1)}, \theta^{(2)}$ be  the sequences of the points of discontinuity
of these functions, respectively. We shall say that $u_1$ is $\,\epsilon-$ equivalent to $u_2$, and denote
$u_1 \epsilon u_2,$  if  $\theta^{(1)}  \epsilon  \theta^{(2)}$ and
 $|u_1(t) -  u_2(t) | < \epsilon $ for all 
$ t\in R \backslash \cup_i [(\theta_i^{(1)} - \epsilon   , \theta_i^{(1)} +\epsilon) \cup
(\theta_i^{(2)} - \epsilon   , \theta_i^{(2)} +\epsilon)]. $
We also say that $u_1$ belongs to the $\epsilon-$ neighborhood of $u_2,$ and vice versa, denoting  $u_1 \in O(u_2, \epsilon)$ and  $u_2 \in O(u_1, \epsilon),$ respectively.
\end{defn}
\begin{defn} The topology defined on the basis of all $\epsilon-$ neighborhoods of functions from
${\cal PC}_r$ is called $B-$ topology.
It is clear that this  topology is Hausdorff.
\end{defn}
\begin{defn} A number $\tau$ is an $\epsilon-$ translation number of $\phi \in {\cal PC}_r$
if $\phi(t+\tau) \in  O(\phi(t), \epsilon).$
\end{defn}
\begin{defn} A function $\phi \in {\cal PC}_r$ is a Bohr-Wexler almost periodic function, if for arbitrary 
$\epsilon>0$ there exists a relatively dense set of $\epsilon-$ translation numbers of  $\phi.$  If $\theta \in \Theta$ is
a sequence of the moments of discontinuity of $\phi,$  then $\theta$ is a  Wexler sequence.
\end{defn}
We shall denote  by ${\cal BWAP}$  the set of all Bohr-Wexler almost periodic functions.
Let $h = \{h_n\}$ and $T_h \phi \in {\cal PC}_r$ be the limit of the sequence $\phi(t+ h_n), \phi \in {\cal PC}_r,$ in
$B-$ topology, if it exists.
\begin{defn} $\phi \in {\cal PC}_r$ satisfies the Bohner property if every sequence $h^{'}$ contains a subsequence
$h \subseteq h^{'}$ such that there exists $T_h \phi \in {\cal PC}_r.$
\end{defn}
\begin{lem} If $u_1\epsilon_1 u_2,   u_2\epsilon_2 u_3$, then $u_1(\epsilon_1+ \epsilon_2) u_3.$
\label{l6}
\end{lem}
{\it Proof.} Lemma \ref{l1} and  the relations $\theta^{(1)} \epsilon_1 \theta^{(2)},   \theta^{(2)} \epsilon_2 \theta^{(3)}$
imply that   $\theta^{(1)} (\epsilon_1+ \epsilon_2) \theta^{(3)}.$ Moreover, one can easily obtain  that 
\begin{eqnarray}
t \in R \backslash \cup_i [(\theta_i^{(1)}-(\epsilon_1+ \epsilon_2), \theta_i^{(1)}+(\epsilon_1+ \epsilon_2)) \cup
(\theta_i^{(3)}-(\epsilon_1+ \epsilon_2), \theta_i^{(3)}+(\epsilon_1+ \epsilon_2))]
\label{e16}
\end{eqnarray}
implies that $t\in R \backslash \cup_i [(\theta_i^{(j)}-\epsilon_j, \theta_i^{(j)}+\epsilon_j) \cup
(\theta_i^{(2)}-\epsilon_j, \theta_i^{(2)}+\epsilon_j)], j = 1,3.$\\
 That is why  $|u_1(t) - u_3(t)| \leq |u_1(t) - u_2(t)| + |u_2(t) - u_3(t)| < \epsilon_1+ \epsilon_2$ if (\ref{e16}) is valid.
The lemma is proved.
\begin{thm} $\phi \in {\cal BWAP}$ if and only if $\phi$ satisfies the Bohner property. 
  $T_h \phi \in {\cal BWAP}$ for $\phi \in {\cal BWAP},$ if the limit exists.
\end{thm}
{\it Proof.} {\it Necessity.} Assume that $\phi \in {\cal BWAP}, \theta \in \Theta$ is a sequence 
of the points of discontinuity of $\phi,$  and $h'\subset \mathbb R$ is a given sequence.
By Theorem \ref{t1} there exists a subsequence  of $h'$  such that  $T_{h^{'}}\theta = \theta^0$ is a  Wexler sequence,  and without loss of generality 
we assume that it is $h^{'}$
itself.
Consider a sequence  $\epsilon_n$ such that $ \epsilon_n \rightarrow 0, n \rightarrow  \infty,$ and  denote $A_n = \cup_i [\theta_i^0 - \epsilon_n,
\theta_i^0 + \epsilon_n].$ Using the diagonal process  \cite{f}, one can find $h^{(1)} \subseteq h^{'},$
such that $\phi(t+  h_n^{(1)})$ is uniformly convergent to $\phi^{(1)}$ on $A_1$. Then  in the same way we can define a sequence 
$h^{(2)} \subseteq h^{(1)}$ such that  $\phi(t+  h_n^{(2)})$ is uniformly convergent to $\phi^{(2)}$ on $A_2,$ and so on.
Obviously,  $\phi^{(i+1) }= \phi^{(i)}$ on $A_i, A_i \subset A_{i+1}, \cup_i A_i = R\backslash \theta^0.$
Consequently,  $\phi(t+  h_n^{(n)})$ is convergent to $\phi^0 \in {\cal PC}_r$ in $B-$ topology.

Fix $\epsilon >0.$  There exists $n(\epsilon)$ such that for arbitrary $n>n(\epsilon)$ the inequality 
$\phi(t+h_n^{(n)}) \frac{\epsilon}{3} \phi^0(t)$ is valid.  If $\tau$ is  an $ \frac{\epsilon}{3}-$ almost period of
$\phi(t),$ then $\phi(t+h_n^{(n)}+ \tau) \frac{\epsilon}{3} \phi(t+h_n^{(n)})$ and 
$\phi(t+h_n^{(n)}+ \tau) \frac{\epsilon}{3} \phi^0(t+\tau).$ Now, using Lemma \ref{l6}, one can obtain
$\phi^0(t+\tau) \in O(\phi^0, \epsilon).$ 

{\it Sufficiency.} Assume that  $\phi \not \in {\cal BWAP}$. Then, similarly  to the classical case for some $\epsilon_0 >0,$
 we can find a sequence of sections $I_n = [h_n -l_n,h_n +l_n], l_n \geq \max_{m<n} |h_m|, l_1-$ arbitrary, such that
if $\xi \in \cup_n I_n,$ then $\phi(t+\xi) \not \in O(\phi, \epsilon_0).$  Denote $h^{'} = \{h_n\},$ and assume that 
there exists $h \subset h^{'}$ such that $T_{h^{'}} = \phi^0 \in {\cal PC}_r.$   Then there exists $n(\epsilon)$ such  that
 if $m>p> n(\epsilon),$ then 
$\phi(t+h_m) \frac{\epsilon_0}{2} \phi^0(t),$  and $\phi(t+h_p) \frac{\epsilon_0}{2} \phi^0(t).$ By Lemma \ref{l6}, 
$\phi(t+h_m) \epsilon_0 \phi(t+h_p).$ Hence, $\phi(t+(h_m - h_p)) \epsilon_0 \phi(t), $  but $(h_m -h_p) \in I_m.$
The theorem is proved.

Let  sequences $h^{'}, g^{'}$ be  given.  Subsequences $h \subset h^{'}, g \subset  g^{'}$ 
are common subsequences of  $h^{'}, g^{'},$  respectively,  if $h_n= h'_{n(k)},   g_n= g'_{n(k)}$ for some function
$n(k).$
There is analogue  of Theorem 1.17 from \cite{f}, which can be proved similarly to Theorem \ref{t3}.
\begin{thm} $\phi \in {\cal BWAP}$ if and only if for arbitrary    $h^{'}, g^{'}$ there exist common subsequences
$h, g$ such that $T_{h+g}\phi = T_h T_g \phi.$
\end{thm}

\begin{lem} Assume that $f \in {\cal PC}_r,$ and (\ref{e3}) satisfies exponential dichotomy.
Then the system
\begin{eqnarray}
&& x^{'} = A(t)x + f(t)
\label{e17}
\end{eqnarray}
has a unique solution $ x_{0}(t) \in C_{0}(t).$
\label{l8}
\end{lem}
{\it Proof.} Indeed,  similarly to   Lemma \ref{l1} one can check that 
\begin{eqnarray}
&& x_0(t)  = \int_{-\infty}^{\infty} G(t,s) f(s)ds
\label{e18}
\end {eqnarray}
is a solution of (\ref{e17}) and belongs to  $C_0(\mathbb R).$ Assume that $x_1(t)$ is another 
solution of  (\ref{e17}),  bounded on $\mathbb R.$ One can see that
the difference $x_1 - x_0$ is a  continuously differentiable solution of  system (\ref{e3}).
Hence, it is a trivial solution of (\ref{e3}). The lemma is proved.\\
In what follows we assume that
\begin{itemize}
\item[C5)]  $\theta$ is a Wexler sequence.

\end{itemize}
Using the Bohner property one can prove that the following assertion is valid.
\begin{lem} Assume that $\phi(t) \in {\cal AP} (\mathbb R)$, and condition $(C_5)$ is valid. Then $\phi(\beta(t)) \in {\cal BWAP}.$
\label{l9}
\end{lem} 

Assume that $\phi(t) \in {\cal AP}$ and $\psi(t) \in {\cal BWAP}.$

For our convenience, following \cite{amerio}, we shall say that
\begin{itemize} 
\item[(i)]   a sequence $h$ is regular with respect to $\phi(t)$ if  the sequence $\phi(t+h_n)$ is uniformly convergent  on $\mathbb R;$
 \item[(ii)] a sequence $h$ is regular with respect to $\psi(t)$ if the sequence $\psi(t+h_n)$ is convergent in $B-$ topology.
\end{itemize}
Let us denote by ${\cal L}(\phi) ({\cal L}(\psi))$ the set of all sequences regular with respect to 
$\phi (\psi).$
\section{Almost periodic solutions}

\begin{lem} Assume that (\ref{e3}) satisfies exponential dichotomy and $f(t) \in {\cal BWAP}.$
Then there is a unique solution of (\ref{e17}), $x_0(t) \in  {\cal AP(\mathbb R)},$   such that ${\cal L}(x_0)
\subseteq {\cal L}(A,f),$ and $||x_0|| \leq (\frac{K_1}{\sigma_1} +\frac{K_2}{\sigma_2})||f||.$ 
\label{l10}
\end{lem}  
{\it Proof.} By Lemma  \ref{l8}, the function $x_0(t)$  defined by (\ref{e18}) is 
a solution of (\ref{e17}) and $x_0(t) \in C_0(\mathbb R).$ 	It is easy to verify that
$||x_0|| \leq (\frac{K_1}{\sigma_1} +\frac{K_2}{\sigma_2})||f||.$  Since every system in the hull of (\ref{e3})
satisfies exponential dichotomy \cite{cop}, it has a unique bounded solution on $\mathbb R.$
Let sequences $h^{'}$ and  $g^{'}$ be given. There exist common sequences $h \subset h^{'}$ and $g\subset g^{'}$ such that 
$T_{h+g} A = T_{h}T_{g} A, T_{h+g} f= T_{h}T_{g} f,$  and there exist uniform limits on 
compact sets $y = T_{h+g} x_0$ and $z= T_{h}T_{g} x_0.$  Since $y,z \in C_0(\mathbb R)$ and 
they are solutions of the same equation,  it follows that $y=z.$ By Theorem 1.17 \cite{f},   $x_0(t)$ is  an almost periodic function.
Assume that for a given sequence $h$ we have $T_h A = A^*$ and $T_{h} f = f^*.$ We shall show that
 the limit $T_{h} x_0 $ exists. Indeed,  suppose, on the contrary,  that the limit does not exist.  Then there are two 
subsequences $h^{(1)} \subset h$ and $h^{(2)} \subset h$ such that
$$|| x_0(t+h_{1n}) -  x_0(t+h_{2n})|| \geq \epsilon_0>0.$$
Then  $$|| x_0(t+h_{1n}-h_{2n} ) -  x_0(t)|| \geq \epsilon_0>0,$$
for all $n \in \mathbb N.$
But $T_{h_{1}-h_{2}} A = A,$ and  $T_{h_{1}-h_{2}} f = f.$  Hence, $T_{h_{1}-h_{2}} x_0  = x_0.$
The theorem is proved. 

Using the Bohner property again and Lemma \ref{l9},  one can prove that the following lemma is valid.
\begin{lem} If  $ f \in {\cal AP}(\mathbb R\times G_H^{m}),$  and  $\psi \in {\cal BWAP}, 
\psi: R \rightarrow G_H,$  
then  $ F_{\theta}(\psi(t)) \in {\cal BWAP}$ and  ${\cal L}( F_{\theta}(\psi(t))) \subseteq  {\cal L}(f, \psi).$ 
\label{lem8}
\end{lem}
\begin{thm} Assume that conditions $(C1)-(C5)$ are valid,  and 
$$ \, l m (\frac{K_1}{\sigma_1}+\frac{K_2}{\sigma_2}) <1.$$

Then there exists a unique solution of  (\ref{e1}),
 $\phi(t) \in {\cal AP}(\mathbb R),$  such that   ${\cal L}(\phi) \subseteq {\cal L}(A,f,\beta).$
\label{thm1}
\end{thm}
{\it Proof.} \rm Let $\Psi = \{\psi \in \cal{AP}(\mathbb R)| \cal{L}(\psi) \subseteq \cal{L}(A,F,\beta)\}$
be  a complete metric space with the  $\sup-$norm $||\cdot||_0.$
Define an operator $\Pi$ on $\Psi$ such that 
\begin{eqnarray*}
&&\Pi( \psi(t) )= \int_{-\infty}^{\infty} G(t,s)F_{\theta}(\psi(s))ds.
\label{e12'}
\end{eqnarray*}
Lemma \ref{lem8} implies that $F_{\theta}(\psi(s)) \in \Psi$  and  $\Pi: \Psi \rightarrow \Psi.$
If $\psi_1, \psi_2 \in \Psi,$  then 
\begin{eqnarray*}
&&||\Pi( \psi_1(t)) -   \Pi( \psi_2(t))|| \leq || \int_{-\infty}^{t} X(t)PX^{-1}(s)(F_{\theta}(\psi_1(s)) -
F_{\theta}(\psi_2(s)) ds ||+\\
&&||\int_{t}^{\infty} X(t)(I-P)X^{-1}(s)(F_{\theta}(\psi_1(s)) -
F_{\theta}(\psi_2(s)) ds || \leq \\
&&l m  (\frac{K_1}{\sigma_1}+\frac{K_2}{\sigma_2}) ||\psi_1(t)) -  \psi_2(t)||_0.
\end{eqnarray*}
Thus,  $\Pi$ is a contractive operator and there exists a unique almost periodic
solution of the equation 
\begin{eqnarray*}
&&\psi(t) = \int_{-\infty}^{\infty} G(t,s)F_{\theta}(\psi(s))ds,
\end{eqnarray*}
which is a solution of (\ref{e1}). 
The theorem is proved.
\begin{rem} Lemma \ref{l10}   and Theorem  \ref{thm1} are analogous to the  assertions 
which were obtained in \cite{f} for ordinary differential equations. 
\end{rem}

\section{ Stability}
This  section is concerned with the  problem of  
stability of the almost periodic solution of  system (\ref{e1}). 
We consider a specific initial condition when values of solutions are evaluated only at points from 
sequence $\theta.$ This approach to the stability is natural for 
EPCA \cite {w,w2}. More detailed discussion of the problem can be found in \cite{ap7}.
Denote  by $X(t,s) = X(t)X^{-1}(s)$  the Cauchy matrix of (\ref{e3}).
We will need the following assumptions:
\begin{itemize}
\item[($C_6$)] $\exists \{\sigma, K\} \subset R, K \geq 1, \sigma >0,$ such that 
$ ||X(t,s)|| \leq K \exp ({-\sigma(t-s)}), t\geq s ;$
\item[($C_7$)] $ l < \frac{\sigma}{mK}.$
\end{itemize}
Assume that $p_j  \geq 0, j=\overline{1,m}$ and denote
$\tau = \max\{ \sup_t (t-\theta_{\beta(t)-p_j}), j=\overline{1,m}\} > 0 ,\quad \zeta(l) = 1- \exp (a \tau) K l m (\sigma-a)^{-1},$
where  $a \in R, 0<a<\sigma,$ is fixed.
\begin{itemize}
\item[($C_8$)] $\zeta(l) >0.$
\end{itemize}
Conditions $(C1)-(C7)$ and  Theorem \ref{thm1} imply that there exists a unique solution of (\ref{e1}), $\xi(t) \in {\cal AP}(\mathbb R).$ 

Fix $\epsilon >0 $ and denote $ L(l,\delta) = \frac{K}{\zeta(l)} \delta, $ where $\delta \in R, \delta >0$. 
Take  $\delta$ so small
that $ L(l,\delta)<\epsilon.$
Assume that $t_0 \in \theta.$ Moreover, without any loss of  generality, assume that $t_0 = \theta_0 =0.$  Fix a sequence 
$ \eta^j \in  \mathbb R^n, j = \overline{1,m}, \quad \max||\eta^j|| < \delta.$ 
Denote $p^0 = \max_{\overline{1,m}}p_j,$ and let 
 $\Psi_{\eta}$  be the set of all continuous functions  which are defined on  $[\theta_{-p^0}, \infty).$ And  if  $ \psi \in \Psi_{\eta}$  
then: $1) \, \psi(\theta_{-p^j}) = \eta^j, j=\overline{1,m}; \quad 2)\, \psi(t)$ is  
uniformly continuous  on $[0,+ \infty);$
 and  
3) \, $ || \psi(t)|| \leq L(l,\delta)
 \exp(-a t)$
 if $ t\geq 0.$

Consider the following EPCAG and the initial condition

\begin{eqnarray}
&& \frac{dv}{dt} = A(t)v + w(t,v(\theta_{\beta(t) - p_1}) ,v(\theta_{\beta(t) - p_2}), \ldots, v(\theta_{\beta(t) - p_m})),\nonumber\\
&& v(s) = \eta^j, j=\overline{1,m},
\label{e7}
\end{eqnarray}
where  
\begin{eqnarray*}
&&w(t,v(\theta_{\beta(t) - p_1}) ,v(\theta_{\beta(t) - p_2}), \ldots, v(\theta_{\beta(t) - p_m})) = \\
&&f(t,
\xi(\theta_{\beta(t) - p_1})+
v(\theta_{\beta(t) - p_1}) ,
\xi(\theta_{\beta(t) - p_2})+\\
&& v(\theta_{\beta(t) - p_2}), \ldots, \xi(\theta_{\beta(t) - p_m})+
v(\theta_{\beta(t) - p_m})) -\\
&&f(t,
\xi(\theta_{\beta(t) - p_1}),
\xi(\theta_{\beta(t) - p_2}), \ldots, \xi(\theta_{\beta(t) - p_m})),
\end{eqnarray*} 
and 
$w $ satisfies $w(t,0) = 0,$
\begin{eqnarray*}
&& ||w(t,v_1) - w(t,v_2)|| \leq l \sum_{j=1}^{m} ||v_1^j - v_2^j ||,
\end{eqnarray*}
 $v_i = (v_i^1,\ldots,v_i^{m}) \in R^{nm}, i=1,2.$
 The following definition is an adapted version of a definition 
 from \cite{cw1}
 \begin{defn}  A function $v(t)$ is a solution of the initial value problem (\ref{e7})
on the interval $[\theta_{-p^0}, \infty)$   if
the following conditions are fulfilled:
\begin{enumerate}
\item[(i)] $v(\theta_{-p_j}) = \eta^j, j=\overline{1,m};$
\item[(ii)] $v(t)$ is continuous on $[\theta_{-p^0}, \infty);$
\item[(iii)] the derivative $v'(t)$ exists at each point 
$t \in [0, \infty)$ with the possible exception
of the points $\theta_j, j  \ge 0,$ where one-sided derivatives exist;
\item[(iv)] equation (\ref{e7}) is satisfied  by $v(t)$ on each interval 
 $[\theta_j,\theta_{j+1}),  j  \ge 0.$ 

\end{enumerate}
\label{defn1}
\end{defn}
It is obvious that 
$v'(t)$ is  the restriction on 
$[0, \infty)$ of a function from ${\cal PC}_r$  and the  last definition can be used for EPCAG (\ref{e1}),   too.

Denote $\bar \theta = \sup_{i}(\theta_{i+1} - \theta_i).$ There exists a positive number $M$ 
such that $||X(t,s)|| \le M$ if $|t-s| \le \bar \theta.$

Assume additionally that
	
\begin{itemize}
	\item [$(C_9)$] $\quad M \bar \theta m l < 1.$
\end{itemize}
\begin{thm}
Assume that  $(C_1)-(C_3), (C_5)-(C_9)$ are valid. Then there exists a unique solution of  the initial value problem (\ref{e7}), $v(t)\in \Psi_{\eta}.$
\label{t8}
\end{thm}
{\it Proof.} \rm  
Similarly to Lemma \ref{l1} we can check that the initial value problem   is equivalent to the following integral equation 
\begin{eqnarray}
&&v(t) = X(t,0)\eta^0+ \int_0^t X(t,s)F_w(v(s))ds,\nonumber\\
&&v(\theta_{-p_j}) = \eta^j, j=\overline{1,m},
\label{ta-ta}
\end{eqnarray}
where $F_w(v(s)) = w(s,v(\theta_{\beta(s) - p_1}) ,v(\theta_{\beta(s) - p_2}), \ldots, v(\theta_{\beta(s) - p_m})).$
 Define on $\Psi_{\eta}$ an operator $\Pi$ such that if $\psi \in \Psi_{\eta},$ then 
\begin{eqnarray*}
&&  \Pi \psi= \left\{\begin{array}{rr} \psi, t \in [\theta_{-p^j},0],\\ \nonumber\\
 X(t,0)\eta^0 + \int_0^t X(t,s)F_{\omega}(\psi(s))ds,  t \ge 0\, .
\end{array}\right.\\
\end{eqnarray*}
We shall show that $\Pi : \Psi_{\eta} \rightarrow \Psi_{\eta}.$   
Indeed, for $ t\geq 0$ it is true that 
$$||\Pi\psi|| \leq K\exp(-\sigma t) \delta +  \int_{0}^t K\exp(-\sigma (t-s)) l L(l,\delta) \sum_{j=0}^{m} \exp(-a \theta_{\beta(s)- p_j})ds
\leq $$ $$\exp(-a t) [K\delta + \frac{m \exp(a \tau) K l L(l,\delta)}{\sigma-a}] = L(l,\delta) \exp (-a t).$$
Differentiating $\Pi\psi$ on $[0,\infty),$ it is easy to show that $[\Pi \psi]'$ exists on $[0,\infty)$ except possibly  on a countable
set of isolated points of discontinuity of the first kind,  and that it is 
bounded on  $[0,\infty).$  Hence, $\Pi\psi$ is a uniformly continuous on  $[0,\infty)$ function. \\
Let $\psi_1, \psi_2 \in \Psi_{\eta}.$
Then 
$$||\Pi\psi_1 -  \Pi\psi_2|| \leq  \int_{0}^t b\exp(-a(t-s)) l m ||\psi_1 - \psi_2||_{1} ds \leq 
\frac{Klm}{\sigma} \sup_{t\geq 0}||\psi_ 1- \psi_2||.$$
Using a contraction mapping argument,  one can conclude that there exists 
a unique fixed point $v(t, \eta)$ of the operator $\Pi :\  \Psi_{\eta} \rightarrow \Psi_{\eta}$
which is a solution   of (\ref{e7}). To complete the proof we should show that there exists no  solution of the problem out of  $\Psi_{\eta}.$ 

Consider first the interval $[\theta_0, \theta_1].$ Assume that on the interval (\ref{e7})
has two different solutions $v_1,v_2$ of the problem. Obviously,  their difference $w = v_1-v_2$ is again a solution of the equation.
Denote $\bar m = \max_{[\theta_0, \theta_1]}||w(t)||,$ and assume, on contrary, that $m >0.$ 
We have that on the interval 
	\[ ||w(t)|| =  ||\int_0^t X(t,s)F_w(w(s))ds || \le M l \bar \theta m \bar m.
\]
 The last inequality contradicts condition $(C_9).$ Now, using  induction, one can easily prove the
 uniqueness for all $t \ge 0.$
The theorem is proved.

Denote $\phi = \{\phi^j\}, j = \overline{1,m},$ a sequence of vectors from $\mathbb R^n.$
Let $x(t, \phi)$ be a solution
of (\ref{e1})  such that   $x(\theta_{-p^j}, \phi) = \phi^j,  j = \overline{1,m}.$  

\begin{defn} \rm   The almost periodic solution   $\xi(t)$ of (\ref{e1}) is said to be  exponentially stable if there exists a number 
$a \in R, a >0,$ such that for every   $\epsilon > 0$ there exists a number 
$\delta = \delta (\epsilon),$ such that 
 the inequality $\max_{j = \overline{1,m}}
 ||\phi^j - \xi(\theta_{-p^j})|| < \delta$  implies 
 $||x(t, \phi)-\xi(t)|| < \epsilon \exp(-at))$ for all $t \geq 0.$
\label{def4}
\end{defn}

Consider now a solution  $x(t, \phi)$ of (\ref{e1}) such that  $\max_{j = \overline{1,m}}
 ||\phi^j - \xi(\theta_{-p^j})|| < \delta.$ 
  Since   the solution  $v(t)$ of the equation (\ref{e7}),
 satisfying $v(\theta_{-p^j}) = \phi^j - \xi(\theta_{-p^j}),  j = \overline{1,m},$  exists,  $x(t, \phi)= \xi(t) + v(t), t \in [\theta_{-p^0}, \infty),$ and $x(t, \phi)$ is uniquely continuable to $\infty.$
 Thus, the following theorem is proved.
\begin{thm}
Assume that  $(C_1)-(C_3), (C_5)-(C_8)$ are valid. Then the almost periodic solution $\xi(t)$ of (\ref{e1}) is  exponentially stable. 
\end{thm}

\section{The Example: positive almost periodic  solutions of the logistic differential equation}
\label{positive}
In \cite{s} G. Seifert applied  the reduction to  discrete equations to the following logistic equation $ \dot{N(t)}= N(t)(a(t) - f(N([t]))), t > 0,$  with positive coefficient $a(t).$  In his paper the conditions on the equation  which guarantee the existence of an
asymptotically stable almost periodic solution were found, and it was proposed to solve a
similar problem for an equation of the following type
\begin{eqnarray}
&& N'(t)= N(t)(a(t) - \sum_{j=0}^{n-1}f_j(N([t-j]))), n > 1.
\label{seif1}
\end{eqnarray}
In what follows we  propose a particular solution  of the problem. 

Consider the problem of existence of positive almost periodic solutions of  the following equation
\begin{eqnarray}
&&x'(t)=  x(t)[a(t)-f(x([t - p_1]), x([t- p_2]), \ldots, x([t - p_m]))ds, 
\label{k2}
\end{eqnarray}
where $x \in \mathbb R,$ integers $p_j, j = \overline{1,m},$ are fixed,
$a(t) \in {\cal AP}(\mathbb R), a(t) >0, t \in \mathbb R.$  We  also assume that
$f(0) = 0, \, f$ is positive for positive values of arguments and it is Lipshitzian with positive constant $l$ in all arguments.It is known that  $\theta_i = i, i \in \mathbb Z,$ is a Wexler sequence \cite{hw, sp}. One can easily check that (\ref{seif1}) is of type (\ref{k2}). 

Denote $$M(a) = \lim_{T \rightarrow \infty}\frac{1}{T} \int_{-T}^T
a(s)ds.$$
Assume that $M(a) >0.$ Using the Strengthened Mean Value Theorem \cite{bohr}, 
one can easily verify that there exist positive constants $K, \sigma,$
such that
\begin{eqnarray}
&&  \exp(\int^t_s a(u)du)\leq K \exp(\sigma(t-s)), t \leq s.
\end{eqnarray}
 
By Lemma \ref{l1} the equivalent integral equation is 
\begin{eqnarray}
x(t)=  \int_t^{\infty} \exp(\int^t_s a(u)du) x(s)f(x([s - p_1]), x([s- p_2]), \ldots, x([s - p_m]))ds. 
\label{k6}
\end{eqnarray}

Fix $H \in \mathbb R, 0<H,$ and denote $G = \{(z_0,z_1,z_2, \ldots, z_m) \in \mathbb R^{k+1}\mid z_j \leq H, j = \overline{0,m}\}.$   
Denote $\mu = \sup_{G} z_0f(z_1,z_2, \ldots, z_m).$
We assume that $0<\mu.$  Consider the set of functions 
$\Psi = \{ \psi \in {\cal{AP}(\mathbb R)}\mid 0 \le \psi(t)\leq H , t \in \mathbb R. \}$
 
Define an operator $\Pi$ on the set $ \Psi$ such that if $\psi \in \Psi$ then 
\begin{eqnarray*}
&& \Pi \psi=  \int_t^{\infty} \exp(\int^t_s a(u)du)\psi(s) f(\psi ([s - p_1]), \psi  ([s- p_2]), \ldots, \psi([s - p_k]))ds. 
\end{eqnarray*}

We have  that $ \Pi \psi(t)\leq \frac{K \mu}{\sigma}.$
If we assume that $\frac{K \mu}{\sigma}\leq H$ and $ \frac{K}{\sigma}[lH + \mu] < 1,$
then in the same way as  in Theorem
\ref{thm1} existence of a solution $\psi^0(t)\in \Psi$ can be proved.


\begin{thebibliography}{99} 
 
\bibitem{aw1} A. R. Aftabizadeh, J. Wiener and J.- M. Xu, {\it Oscillatory and periodic solutions of delay differential equations
with piecewise constant argument,}  Proc. Amer. Math. Soc.  99 (1987), 673-679.
\bibitem{amerio} L. Amerio and G. Prouse, \lq\lq Almost-periodic functions and functional equations,\rq\rq \,  Van Nostrand Reinhold Company, new York, 1971.
\bibitem{ap} M. U. Akhmetov and N. A. Perestyuk, {\it Periodic and almost-periodic solutions of strongly nonlinear impulse systems,}   J. Appl. Math. Mech., 56  (1992),  829--837. 
\bibitem{ap1} M. U. Akhmetov and N. A. Perestyuk, {\it
Almost-periodic solutions of nonlinear impulse systems,}   Ukrainian Math. J., 41 (1989), 259-263.
\bibitem{aps} M. U. Akhmetov, N. A. Perestyuk and A.M. Samoilenko,  \lq\lq Almost-periodic solutions of differential equations with impulse action,\rq\rq $\;$ (Russian) Akad. Nauk Ukrain. SSR Inst.,
Mat. Preprint, 1983, no. 26, 49 pp.
\bibitem{ap6} M. U. Akhmet, {\it Existence and stability of almost-periodic  solutions of quasi-linear differential equations with deviating argument,}  Applied mathematics Letters,  17 (2004), 1177-1181.  \bibitem{ap7}M. U. Akhmet,{\it Integral manifolds of  differential equations with piecewise constant
argument of generalized type,} Nonlinear Analysis, Theory, methods and applications, (in press). 
\bibitem{aa} A. Alonso,  and J.  Hong,  {\it Ergodic type solutions of
differential equations with piecewise constant arguments,}   Int. J. Math. Math. Sci.,  28 (2001), 609-619.
\bibitem{aa1} A. Alonso, J.  Hong and R. Obaya,  {\it Almost-periodic type
solutions of differential equations with piecewise constant argument via almost periodic type sequences,}   Appl. Math. Lett.,  13 (2000), 131-137.
\bibitem{bohr}  H. Bohr, \lq \lq Almost-periodic functions,\rq \rq  Chelsea Publishing Company, New York, 1951.
\bibitem{bus} S. Busenberg and K. Cooke, 
\lq \lq Vertically transmitted diseases, Models and dynamics. Biomathematics, 23,\rq\rq Springer-Verlag, Berlin, 1993.
\bibitem{b} T.A. Burton and Tetsuo Furumochi, {\it Fixed points and problems in stability
theory for ordinary and functional differential equations,} Dynamic Systems and Appl., 10  (2001), 89-116.
\bibitem{cw1} K. L. Cooke and J. Wiener, {\it Retarded differential equations with piecewisw constant delays,} J. Math. Anal. Appl.,  99  (1984), 265-297.
\bibitem{cord1} C. Corduneanu, Almost periodic functions, 
Interscience Publishers, New York, 1961.
\bibitem{cop} W.A. Coppel, \lq\lq Dichotomies in stability theory, \rq\rq $\;$  Lecture notes in
mathematics, Springer-Verlag, Berlin, Heidelberg, New York, 1978.
\bibitem{dn} S. Doss and S.K. Nasr, {\it On the functional equation 
$y^{'} = f(x, y(x), y(x+h)), h>0,$}   Amer. Journ. Math., 75  (1953),713 - 716.  
\bibitem{hw} A. Halanay and D. Wexler 
 \lq\lq Qualitative theory of impulsive systems,\rq\rq $\;$ 
(Romanian) Republici Socialiste Romania, Bucuresti, 1968. 
\bibitem{h2}  J. Hale, \lq\lq Functional differential equations,\rq\rq $\;$    Springer-Verlag,
New York, Heidelberg, Berlin, 1971.
\bibitem{g} K. Gopalsamy, \lq\lq Stability and oscillations in delay differential equations,\rq\rq $\;$ Kluwer Academic Publishers Group, Dordrecht, 1992.
\bibitem{fil}   A. F. Filippov,  \lq\lq Differential equations with discontinuous righthand
   sides,\rq\rq $\;$ Mathematics and its Applications (Soviet Series), 18. Kluwer
   Academic Publishers Group, Dordrecht, 1988.
\bibitem{f} A.M. Fink, \lq\lq Almost-periodic differential quations, \rq\rq $\;$  Lecture notes in mathematics, Springer-Verlag, Berlin, Heidelberg, New York, 1974.
\bibitem{e} L.E. El'sgol'ts,  \lq\lq Introduction to the theory of differential equations with deviating arguments, \rq\rq $\;$ 
Holden-Day, Inc, San Francisco, London, Amsterdam, 1966.
\bibitem{kol} A.N. Kolmogorov, {\it On the Skorohod convergence}  (Russian. English summary),  Teor. Veroyatnost. i Primenen.,  1 (1956), 239-247.
\bibitem{kras}   M. A. Krasnosel'skii, V. Sh. Burd and Yu. S. Kolesov,  \lq\lq Nonlinear almost periodic oscillations,\rq\rq John Wiley Sons, New York-Toronto, 1973.
\bibitem{ky} T. K\"upper and R. Yuan,  {\it On quasi-periodic solutions of
   differential equations with piecewise constant argument,} J. Math. Anal. Appl.,  267  (2002),  173-193.
\bibitem{K} Yu. A. Kuznetsov, \lq\lq Elements of Applied Bifurcation Theory,\rq\rq Springer-Verlag, New-York, Berlin, Heidelberg, 1995.  
\bibitem{m} R.M. May and G.F. Oster,  {\it  Bifurcations and dynamic complexity in simple ecological models,} Amer. Natural,  110 (1976) 573-599.
\bibitem{p} G. Papaschinopoulos, {\it Some results concerning a class of differential equations
with piecewise constant argument,}  Math. Nachr.  166 (1994), 193-206.
\bibitem{p1} G. Papaschinopoulos,  {\it Linearisation near the integral manifold for a system of 
differential equations with piecewise constant argument,}  J. of Math. Anal. and Appl.
 (1997), 317-333.
\bibitem{shah} S.M. Shah and J. Wiener, {\it Advanced differential equations with piecewise constant argument deviations,}  Internat. J. Math. Sci,  6 (1983), 671-703. 
 \bibitem{sp} A. M. Samoilenko and N. A. Perestyuk,
 \lq\lq Impulsive Differential Equations,\rq\rq $\;$  World
Scientific, Singapore, 1995.
\bibitem{s} G. Seifert,  {\it Almost periodic solutions of certain  differential equations with piecewise constant delays and almost periodic time dependence,}  J. Differential equations,  164 (2000), 451-458.
 \bibitem{skor}  A.V. Skorokhod, {\it  Limit theorems for random processes,}
Theory Probab. Appl., 39 (1994),289-319,
\bibitem{slyus} V.E. Slyusarchuk, {\it Bounded solutions of impulsive systems,}  Differentsial'nye Uravneniya,  19 (1983) 588-596.
\bibitem{w} D. Wexler, {\it Solutions p\'eriodiques et
presque-p\'eriodiques des syst\'emes d'\'equations diff\'eretielles  lin\'eaires en distributions,} J. Differential Equations., 2 (1966), 12-32.
\bibitem{w2} J. Wiener, {\it  Generalized solutions of functional differential equations,} World Scientific, Singapore (1993).
\bibitem{wl} J. Wiener and  V. Lakshmikantham, {\,\it A damped oscillator with
   piecewise constant time delay,}  Nonlinear Stud.,  7  (2000), 78--84.            
\bibitem{y1} Yuan Rong, {\it On the spectrum of almost periodic solution of second order scalar functional differential equations with piecewise constant argument,}  J. Math. Anal. Appl., 303 (2005) 103--118.  
\bibitem{ym} Muroya, Yoshiaki,  {\it Persistence, contractivity and global stability in logistic equations with piecewise constant delays,}  J. Math. Anal. Appl.,   270 (2002), 602-635. 
\bibitem{z} C. Zhang, \lq\lq Almost periodic type functions and ergodicity, 
Science press,\rq\rq \\ Beijing/New York, Kluwer Academic Publishers, Dordecht/ Boston/ London, 2003.

\end{thebibliography}
\end{document}